\documentclass[thmsa]{article}
\usepackage{amssymb}
\usepackage{amsmath,amsthm,amsfonts}  
\usepackage[applemac]{inputenc}  
\usepackage{graphicx}
\usepackage{amssymb}
\usepackage{amsmath}
\usepackage{amsthm}
\usepackage{color}

\newtheorem{prop}{Proposition}[section]

\newcommand{\Proof}{\noindent{\bf Proof.\,\,\,}}
\renewcommand{\qed}{\hfill $\square$}

\newcommand{\ket}[1]{|\kern.3ex#1\kern.3ex\rangle}
\newcommand{\bra}[1]{\langle\kern.3ex #1 \kern.3ex|}
\newcommand{\scalar}[2]{\langle\kern.3ex #1 \kern.3ex|\kern.3ex#2\kern.3ex\rangle}

\def\N{\mathbb{N}}

\def\Z {\mathbb{Z}}

\definecolor{hervecolor}{rgb}{0.8,0,0.7}

\begin{document}

\author{
E.M.F.  Curado$^{\mathrm{a,b,c}}$,  J. P. Gazeau$^{\mathrm{b}}$, 
 Ligia M.C.S. Rodrigues$^{\mathrm{a}}$(\footnote{e-mail:
 evaldo@cbpf.br,
gazeau@apc.univ-paris7.fr, ligia@cbpf.br} )\\
\emph{$^{\mathrm{a}}$ Centro Brasileiro de Pesquisas Fisicas,}\\
\emph{Rua Xavier Sigaud 150, 22290-180 - Rio de Janeiro, RJ, Brazil}\\
\emph{$^{\mathrm{b}}$ Laboratoire APC,
Univ Paris  Diderot, Sorbonne Paris Cit\'e, }   
\emph{75205 Paris, France} \\
\emph{$^{\mathrm{c}}$ Instituto Nacional de Ciencia e Tecnologia - Sistemas Complexos}}

\title{A note about combinatorial sequences and Incomplete Gamma function}
\author{
H. Bergeron$^{\mathrm{a}}$, 
E.M.F.  Curado$^{\mathrm{b,c}}$,  \\
J. P. Gazeau$^{\mathrm{b,d}}$, 
 Ligia M.C.S. Rodrigues$^{\mathrm{b}}$(\footnote{e-mail: herve.bergeron@u-psud.fr, 
 evaldo@cbpf.br,
gazeau@apc.univ-paris7.fr, ligia@cbpf.br} )\\
\emph{$^{\mathrm{a}}$
Univ Paris-Sud, ISMO, UMR 8214, 91405 Orsay, France} \\
\emph{   $^{\mathrm{b}}$ Centro Brasileiro de Pesquisas Fisicas   } \\
\emph{   $^{\mathrm{c}}$ Instituto Nacional de Ci\^encia e Tecnologia - Sistemas Complexos}\\
\emph{  Rua Xavier Sigaud 150, 22290-180 - Rio de Janeiro, RJ, Brazil  } \\
\emph{  $^{\mathrm{d}}$ APC, UMR 7164,}\\
\emph{ Univ Paris  Diderot, Sorbonne Paris Cit\'e,  }  
\emph{75205 Paris, France} 
}

\maketitle
\abstract{ In this short note we present a set of interesting and useful properties of a one-parameter family of sequences including factorial and  subfactorial, and their relations to the Gamma function and the incomplete Gamma function. }

\tableofcontents

\section{Introduction}
\label{intro}
The purpose of this note is to show a set of properties of a  one-parameter family of sequences $\mathcal{T}_n( z)$, $z$ is a complex variable, which are apparently not known in the literature. Among these sequences we have factorials and subfactorials and they are seen to be related to the incomplete Gamma function. This family was introduced in a recent paper about generalized binomial distribution \cite{bergeronetal2013}.
\section{The function $T_n(x,y)$}
\label{subfctnbs}

Let us consider the following two-variable polynomial:
\begin{equation}
\label{polTxy}
T_n(x,y) := \sum_{m=0}^n \binom{n}{m} (m+x)^m\,(n-m+y)^{n-m}\, , 
\end{equation}
which is clearly symmetrical in $x$, $y$. Actually, we have the stronger result:
\begin{prop}
\label{nodepx}
The polynomial $T_n(x,y)$ depends on the sum $z=x+y$ only:
\begin{equation}
\label{polTnTn}
T_n(x,z-x)= \sum_{m=0}^n \binom{n}{m} (m+x)^m\,(n-m+z-x)^{n-m}\equiv \mathcal{T}_n( z)\, ,
\end{equation}
where $x$ can be arbitrarily chosen with the convention that $0^0\equiv 1$ in the cases $x=0$ or $x=z$.
\end{prop}

\Proof
Let us prove by recurrence on $n$ that $T_n(x,z-x)$ does not depend explicitly on $x$ for all $z$. For that, we take the partial derivative of $T_n(x,z-x)$ with respect to $x$. We have successively, by appropriately shift the summation variable, 
\begin{align*}
\frac{\partial}{\partial x} T_n(x,z-x)=& \sum_{m=0}^n \binom{n}{m}
\left\lbrack m (m+x)^{m-1}\,(n-m+z-x)^{n-m} \right.\\&\left. - (n-m) (m+x)^{m}\,(n-m+z-x)^{n-m-1}\right\rbrack\\
=& n\, \sum_{m=0}^{n-1} \binom{n-1}{m} (m +1 +x)^{m}\,(n-1-m+z-x)^{n-1-m} \\
& - n\, \sum_{m=0}^{n-1} \binom{n-1}{m} (m + x)^{m}\,(n -1-m+z+1-x)^{n-1 -m } \\
=& n \left\lbrack T_{n-1}(1+x,z -x) - T_{n-1}(x,z+1-x)\right\rbrack\,. 
\end{align*}
Now, we note that $T_0(x, z-x)=1$ and so does not depend explicitly on $x$. Suppose that the property ``$T_k(x,z-x)$ does not depend explicitly on $x$ for all $z$'' holds true  for all integer $1\leq k\leq n-1$. Then $T_{n-1}(1+x,z -x) = T_{n-1}(x,z+1-x)$ and this implies  that $\partial/\partial x T_n(x,z-x) = 0$, i.e., $T_n(x,z-x)$ does not depend explicitly on $x$ for all $z$.
\qed

\section{The function $\mathcal{T}_n(z)$}
\label{tnz}

Thanks to Proposition \ref{nodepx}, it is now possible to find a convenient  alternate expression of $T_n(x,z-x)\equiv \mathcal{T}_n(z)$. First, let us establish a recurrence formula.

\begin{prop}
\label{recTn}
The polynomial $\mathcal{T}_n(z)$ satisfies the following recurrence formula:
\begin{equation}
\label{recTnz}
\mathcal{T}_n(z) = (z+n)^n + n\, \mathcal{T}_{n-1}(z+1) \, , \quad \mathcal{T}_0(z)=1\,.
\end{equation}
\end{prop}

\Proof
Taking $x=0$ in  \eqref{polTnTn},  let us write that expression as
\begin{align*}
\mathcal{T}_n( z) &= \sum_{m=0}^n \binom{n}{m} (m)^m\,(n-m+z)^{n-m}\\
&= (n+z)^n + n\sum_{m=1}^{n}\binom{n-1}{m-1}\times\\
&\times(m-1+1)^{m-1}(n-1-(m-1)+z +1 -1)^{n-1-(m-1)}\\
&= (n+z)^n + n T_{n-1}(1,z+1-1) = (n+z)^n +n \mathcal{T} _{n-1}(z+1)\,  ,
\end{align*}
where we have applied Proposition \ref{nodepx} with $x=1$. 
\qed

It is then straightforward to  deduce from this formula the following result. 

\begin{prop}
\label{alterTn1}
The polynomial $\mathcal{T}_n(z)$ admits the following expansion in powers of $ (z+ n)$:
\begin{equation}
\label{z+n}
\mathcal{T}_n(z) = \sum_{k=0}^n \frac{n!}{k!} (z+ n)^k\, .
\end{equation}
\end{prop}
We note in particular the interesting corollary proving that the family $\left(\mathcal{T}_n(z)\right)_{z\in \Z}$ of integer sequences includes the factorial:
\begin{equation}
\label{z=-n}
\mathcal{T}_n(-n) = n!\, ,
\end{equation}
\label{gamma}
and consequently also the Gamma function as 
\begin{equation}
\label{gamma1}
\mathcal{T}_n(-n) = \Gamma(n+1)\,  .
\end{equation}

The expression \eqref{z+n} allows to easily derive two asymptotic behaviors of $\mathcal{T}_n(z)$: 
\begin{align}
\label{larzTn}
\mbox{At large}\ z    & \quad \mathcal{T}_n(z) \sim z^n \, ,\\
\label{larzTn} \mbox{At large}\ n  &\quad   \mathcal{T}_n(z) \sim n^n\,.
\end{align}
Indeed, in the latter the dominant term of the sum is for $k=n$ (the term $k=0$ is not dominant 
due to the decreasing exponential factor in the Stirling's formula 
$n!\sim n^n \sqrt{2\pi n} \, e^{-n}$). 
This indicates that the dominant term of the sum  in \eqref{z+n} is $(z+n)^n$. 

A second expression of $\mathcal{T}_n(z)$ is also of interest. 
\begin{prop}
\label{alterTn2}
The polynomial $\mathcal{T}_n(z)$ admits the following expansion in powers of $ (z+ n +1)$:
\begin{equation}
\label{z+n+1}
\mathcal{T}_n(z) = \sum_{k=0}^n a_k (z+ n +1)^k\, , 
\end{equation}
where the coefficients $a_k$ are given by
 \begin{equation}
\label{akTn}
a_k = \binom{n}{k}\,\mathcal{T}_{n-k}(k-n-1)\equiv\binom{n}{k}d_{n-k}\, . 
\end{equation}
\end{prop}

\Proof
Applying formula $$a_k= \frac{1}{k!}\frac{d^k}{dz^k} \left.\mathcal{T}_n(z) \right\vert_{z=-n-1}$$ to the original expression \eqref{polTnTn} where we put $x=0$, we  easily derive
\begin{align*}
\frac{1}{k!}\frac{d^k}{dz^k} \left.\mathcal{T}_n(z) \right\vert_{z=-n-1}& = \sum_{m=0}^{n-k}(-1)^{n-m-k}\binom{n}{m}\binom{n-m}{k}m^m(1+m)^{n-m-k}\\
&= \binom{n}{k}\sum_{m=0}^{n-k}(-1)^{n-m-k}\binom{n-k}{m}m^m(1+m)^{n-m-k}\\
&= \binom{n}{k}\,\mathcal{T}_{n-k}(k-n-1)\, .
\end{align*}
\qed

The sequence of numbers $\left(d_n \equiv \mathcal{T}_{n}(-n-1)\right)_{n\in \N}$, for which \eqref{z+n} gives
\begin{equation}
\label{ }
d_n =  \sum_{k=0}^n \frac{n!}{k!} (-1)^k\, , 
\end{equation}
and whose first terms are 1, 0, 1, 2, 9, 44, ..., is well known for more than 3 centuries \cite{montmort1713,euler1809,comtet74,desarmenien82,oeis}. Its OEIS name is A000166. Its numbers are named   \textit{subfactorial} (and then denoted as $!n$) or \textit{rencontres numbers}, or \textit{derangements}, since $d_n$ is the number of permutations of $n$ elements with no fixed points. They obey the recurrence relations (Euler) $d_n= (n-1)d_{n-1} + (-1)^n$ and $d_n= n(d_{n-1} + d_{n-2})$. Their generating function is 
\begin{equation*}
D(x) = \sum_{n=0}^{\infty}d_n\,\frac{x^n}{n!}= \frac{e^{-x}}{1-x}\,,
\end{equation*}
and their asymptotic behavior at large $n$ is 
\begin{equation*}
\underset{n\to \infty}{\lim} \frac{d_n}{n!}= \frac{1}{e}\,. 
\end{equation*}

\section{The relation to the incomplete Gamma function}
\label{incgam}

Let us finally establish the relation between the polynomial $\mathcal{T}_n(z) $ and the incomplete Gamma function \cite{abramovitz}.

\begin{prop}
\label{Tn2incGfct}
The polynomial $\mathcal{T}_n(z)$ is related to the incomplete Gamma function 
\begin{equation}
\label{incgamfct}
\Gamma(a,x) = \int_{x}^{\infty}t^{a-1}\,e^{-t}\,dt\, , \quad \mathrm{Re}(a) >0\, , 
\end{equation}
as follows:
\begin{equation}
\label{igamnz}
\mathcal{T}_n(z) = e^{z+n}\, \Gamma(n+1,z+n) \,.
\end{equation}
\end{prop}
Note that in the particular case that $z = -n$, from (\ref{igamnz}) we reobtain (\ref{gamma1}).
\proof
It is straightforward to prove, by integration by part, the recurrence formula obeyed by the incomplete Gamma function:
\begin{equation*}
\label{recincgamfct}
\Gamma(a,x)= e^{-x}\,x^{a-1} + (a-1) \Gamma(a-1,x)\Leftrightarrow e^{x}\Gamma(a,x)= \,x^{a-1} + (a-1) e^x\Gamma(a-1,x)\, .
\end{equation*}
Applied to polynomial $\mathcal{T}_n(z)$ this formula gives 
\begin{align*}
\mathcal{T}_n(z) &= e^{z+n}\, \Gamma(n+1,z+n)=  (z+n)^n + n e^{z+n}\, \Gamma(n,z+1+n-1) \\
&= (z+n)^n + n  \mathcal{T}_{n-1}(z+1)\, , 
\end{align*}
and this is precisely \eqref{recTnz} with the same initial condition $\mathcal{T}_0(z)= e^{z}\, \Gamma(1,z)=1$. 
\qed

\section{Comment}
\label{comment}
It is probable that $\mathcal{T}_n(z) $ has a combinatorial interpretation (and alternate simpler expression) for each $z\in \Z$. For instance the elements of the sequence $\left(\mathcal{T}_{n}(1)\right)_{n\in \N}$, whose first terms are 1, 3, 17, 142, 1569,... are the numbers of connected functions on $n$ labeled nodes, and we also have (OEIS number: A001865, see \cite{oeis} for references)
\begin{equation*}
\mathcal{T}_{n}(1) = \sum_{k=0}^n \frac{n!}{k!}(n+1)^k =
e^{n+1}\int_{n+1}^{\infty} x^n e^{-x} \, dx\, . 
\end{equation*}

We can present other examples. For $z=2$: the elements of the sequence $\left(\mathcal{T}_{n}(2)\right)_{n\in \N}$, whose first terms are 1, 4, 26, 236, 2760,..., are the numbers of normalized total height of rooted trees with $n$ nodes (OEIS number: A001863, see \cite{sloane69}). The sequence $(\mathcal{T}_{n}(3))$, whose first terms are 1, 5, 37, 366, 4553,..., is   A129137: $\mathcal{T}_{n}(3)$ is the number of trees on $1,2,3, \dotsc, n \equiv [n]$, rooted at 1, in which 2 is a descendant of 3.  And so on.

\end{document}